\documentclass{article}
\usepackage{amsfonts}
\usepackage{amsmath}

\input{tcilatex}

\begin{document}

\begin{center}
\bigskip

\textbf{ISOTROPIC GEOMETRY OF GRAPH SURFACES ASSOCIATED WITH PRODUCT
PRODUCTION FUNCTIONS IN ECONOMICS}

\bigskip

Muhittin Evren Aydin, Mahmut Ergut

Firat University and Namik Kemal University, Turkey

\bigskip
\end{center}

\textbf{Abstract. }A production function is a mathematical formalization in
economics which denotes the relations between the output generated by a
firm, an industry or an economy and the inputs that have been used in
obtaining it. In this paper, we study the product production functions of 2
variables in terms of the geometry of their associated graph surfaces in the
isotropic $3-$space $\mathbb{I}^{3}$. In particular, we derive several
classification results for the graph surfaces of product production
functions in $\mathbb{I}^{3}$ with constant curvature.

\bigskip



\section{Introduction}

In economics, a \textit{production function} is a non-constant positive
function that species the output of a firm, an industry, or an entire
economy for all combinations of inputs. Explicitly, it is a map of class $%
C^{\infty }$ which has non-vanishing first derivatives defined by%
\begin{equation}
\left\{ 
\begin{array}{l}
h:%
\mathbb{R}
_{+}^{n}\longrightarrow 
\mathbb{R}
_{+},\text{\ }\left( x_{1},x_{2},...,x_{n}\right) \longmapsto h\left(
x_{1},x_{2},...,x_{n}\right) , \\ 
\mathbb{R}
_{+}^{n}=\left\{ \left( x_{1},x_{2},...,x_{n}\right)
:x_{i}>0,i=1,...,n\right\} .%
\end{array}%
\right.  \tag{1.1}
\end{equation}%
Here $h$ denotes the quantity of output, $n$ is the number of inputs and the
variables are the inputs (such as capital, labor, raw materials etc.). Some
interesting examples of production functions can be found in \cite{18}.

In order for the production functions to model as well economic reality,
they are required to get some proporties (see e.g. \cite{5,22}). One of the
most important of these proporties is the production function $f$ to be
homogeneous, i.e. there exist a real number $p$ such that%
\begin{equation}
h\left( \lambda x_{1},\lambda x_{2},...,\lambda x_{n}\right) =\lambda
^{p}h\left( x_{1},x_{2},...,x_{n}\right) ,\text{ }\lambda \in 
\mathbb{R}
_{+}.  \tag{1.2}
\end{equation}%
$\left( 1.2\right) $ implies when the inputs are multiplied by same factor,
the output is multiplied by some power of the factor.

If $p<1$ (resp. $p>1$) in $\left( 1.2\right) $, then the production function
is said to have \textit{decreasing }(resp.\textit{\ increasing})\textit{\
return to scale. }If $p=1,$ then it is said to have \textit{constant return
to scale}.

The presence of increasing returns means that a one percent increase in the
usage levels of all inputs would result in a greater than one percent
increase in output; the presence of decreasing returns means that it would
result in a less than one percent increase in output. Constant returns to
scale is the in-between case (cf. \cite{8}).

A.D. Vilcu and G.E. Vilcu \cite{26} completely classified the homogeneous
production functions with constant proportional marginal rate of
substitution. Further, the homogeneous production functions have been
investigated via their associated graph hypersurfaces in \cite{6}-\cite{11}, 
\cite{29}.

The most famous one among homogeneous production functions is \textit{%
Cobb-Douglas production function}, introduced in 1928 by C. W. Cobb and P.
H. Douglas \cite{12}. In original form, it is given as%
\begin{equation*}
Y=bL^{k}C^{1-k},
\end{equation*}%
where $b$ presents the total factor productivity, $Y$ the total production, $%
L$ the labor input and $C$ the capital input.

The generalized \textit{Cobb-Douglas production function }of $n$ variables
is defined by%
\begin{equation*}
h\left( \mathbf{x}\right) =Ax_{1}^{\alpha _{1}}x_{2}^{\alpha
_{2}}...x_{n}^{\alpha _{n}},\text{ }\mathbf{x}=\left(
x_{1},x_{2},...,x_{n}\right) \in 
\mathbb{R}
_{+}^{n}
\end{equation*}%
where $A,\alpha _{1},\alpha _{2},...,\alpha _{n}>0.$ We note that $h$ has
constant return to scale if and only if $\sum_{i=1}^{n}\alpha _{i}=1.$

X. Wang, Y. Fu \cite{28} and A.D. Vilcu, G.E. Vilcu \cite{24,25} classified
the graph hypersurfaces of the generalized Cobb-Douglas production functions
with zero Gauss-Kronocker and mean curvature.

On the other hand, there are some non-homogeneous production functions,
including the famous Spillman-Mitscherlich and transcendental production
functions respectively defined by 
\begin{equation*}
\left\{ 
\begin{array}{l}
h\left( \mathbf{x}\right) =A\left[ 1-\exp \left( a_{1}x_{1}\right) \right]
\cdot \left[ 1-\exp \left( a_{2}x_{2}\right) \right] \cdot ...\cdot \left[
1-\exp \left( a_{n}x_{n}\right) \right] , \\ 
\mathbf{x}=\left( x_{1},x_{2},...,x_{n}\right) \in 
\mathbb{R}
_{+}^{n},\text{ }A>0,\text{ }a_{i}<0,\text{ }i=1,...,n%
\end{array}%
\right.
\end{equation*}%
and%
\begin{equation*}
\left\{ 
\begin{array}{l}
h\left( \mathbf{x}\right) =Ax_{1}^{a_{1}}\exp \left( b_{1}x_{1}\right) \cdot
x_{2}^{a_{2}}\exp \left( b_{2}x_{2}\right) \cdot ...\cdot x_{n}^{a_{n}}\exp
\left( b_{n}x_{n}\right) , \\ 
\mathbf{x}=\left( x_{1},x_{2},...,x_{n}\right) \in 
\mathbb{R}
_{+}^{n},\text{ }A>0,\text{ }a_{i},b_{i}\in 
\mathbb{R}
,\text{ }a_{i}^{2}+b_{i}^{2}\neq 0,\text{ }i=1,...,n.%
\end{array}%
\right.
\end{equation*}%
Such production functions, including the generalized Cobb-Douglas production
function, belong to a more general class of production functions given by%
\begin{equation*}
h\left( \mathbf{x}\right) =\prod\limits_{j=1}^{n}g_{j}\left( x_{j}\right) ,%
\mathbf{x=}\left( x_{1},x_{2},...,x_{n}\right) \in 
\mathbb{R}
_{+}^{n},
\end{equation*}%
where $g_{j}$ is a continuous positive real functions with nonzero first
derivatives, $j=1,...,n.$ H. Alodan et al. \cite{1} called these production
functions \textit{product }production functions. A production function is
said to be \textit{quasi-product}, if it is of the form%
\begin{equation*}
h\left( \mathbf{x}\right) =F\left( \prod\limits_{j=1}^{n}g_{j}\left(
x_{j}\right) \right) ,\mathbf{x=}\left( x_{1},x_{2},...,x_{n}\right) \in 
\mathbb{R}
_{+}^{n},
\end{equation*}%
where $F,g_{j}$ are continuous positive real functions with nonzero first
derivatives, $j=1,...,n.$ A lot of classifications of the quasi-product and
quasi-sum production functions can be found in \cite{1}-\cite{4}, \cite%
{14,27} in terms of the geometry of their graph hypersurfaces.

B.-Y. Chen, et al \cite{7,13} investigated the graph hypersurfaces of the
production models via the isotropic geometry. In the present paper, we
classify the graph surfaces of product production functions of 2 variables
with constant curvature in the isotropic 3-space $\mathbb{I}^{3}.$

\section{Basics on isotropic spaces}

For later use, we provide a brief review of isotropic geometry from \cite{7}%
, \cite{19}-\cite{21}.

The isotropic 3-space $\mathbb{I}^{3}$ is a Cayley--Klein space defined from
a 3-dimensional projective space $P\left( 
\mathbb{R}
^{3}\right) $ with the absolute figure which is an ordered triple $\left(
\omega ,f_{1},f_{2}\right) $, where $\omega $ is a plane in $P\left( 
\mathbb{R}
^{3}\right) $ and $f_{1},f_{2}$ are two complex-conjugate straight lines in $%
\omega $. The homogeneous coordinates in $P\left( 
\mathbb{R}
^{3}\right) $ are introduced in such a way that the absolute plane $\omega $
is given by $X_{0}=0$ and the absolute lines $f_{1},f_{2}$ by $%
X_{0}=X_{1}+iX_{2}=0,$ $X_{0}=X_{1}-iX_{2}=0.$ The intersection point $%
F(0:0:0:1)$ of these two lines is called the absolute point. The group of
motions $G_{6}$ of $\mathbb{I}^{3}$ is a six-parameter group given in the
affine coordinates $x_{1}=\frac{X_{1}}{X_{0}},$ $x_{2}=\frac{X_{2}}{X_{0}},$ 
$x_{3}=\frac{X_{3}}{X_{0}}$ by

\begin{equation}
\left( x_{1},x_{2},x_{3}\right) \longmapsto \left( x_{1}^{\prime
},x_{2}^{\prime },x_{3}^{\prime }\right) :\left\{ 
\begin{array}{l}
x_{1}^{\prime }=a+x_{1}\cos \phi -x_{2}\sin \phi , \\ 
x_{2}^{\prime }=b+x_{1}\sin \phi +x_{2}\cos \phi , \\ 
x_{3}^{\prime }=c+dx_{1}+ex_{2}+x_{3},%
\end{array}%
\right.  \tag{2.1}
\end{equation}%
where $a,b,c,d,e,\phi \in 
\mathbb{R}
.$ Such affine transformations are called \textit{isotropic congruence
transformations }or \textit{i-motions. }It is easily seen from $\left(
2.1\right) $ that i-motions are indeed composed by an Euclidean motion in
the $x_{1}x_{2}-$plane (i.e. translation and rotation) and an affine shear
transformation in $x_{3}-$direction.

In general, the following terminology is used for the isotropic spaces.
Consider the points $\mathbf{x}=\left( x_{1},x_{2},x_{3}\right) $ and $%
\mathbf{y}=\left( y_{1},y_{2},y_{3}\right) .$ The projection in $\mathbf{z}-$%
direction onto $\mathbb{R}^{2},$ $\left( x_{1},x_{2},x_{3}\right)
\longmapsto \left( x_{1},x_{2},0\right) ,$ is called the\textit{\ top view}.
In the sequel, many of metric properties in isotropic geometry (invariants
under $G_{6})$ are Euclidean invariants in the top view such as the
isotropic distance, so-called i-distance. \textit{I-distance }of two points $%
\mathbf{x}$ and $\mathbf{y}$ is defined as the Euclidean distance of their
top views, i.e.,%
\begin{equation}
\left\Vert \mathbf{x}-\mathbf{y}\right\Vert _{i}=\sqrt{\sum_{j=1}^{2}\left(
y_{j}-x_{j}\right) ^{2}}.  \tag{2.2}
\end{equation}%
The i-metric is degenerate along the lines in $\mathbf{z}-$direction, and
such lines are called \textit{isotropic} lines. The plane containing an
isotropic line is called an\textit{\ isotropic plane. }Therefore\textit{, }%
an \textit{isotropic }$3-$\textit{space} $\mathbb{I}^{3}$\ is the product of
the Euclidean $2-$space $\mathbb{R}^{2}$ and an isotropic line\textit{\ }with%
\textit{\ }a degenerate parabolic distance metric.

Let $M^{2}$ be a regular surface immersed in $\mathbb{I}^{3}$ which has no
isotropic tangent planes. Such a surface $M^{2}$ is said to be \textit{%
admissible} and can be parametrized by%
\begin{equation}
r:D\subseteq \mathbb{R}^{2}\longrightarrow \mathbb{I}^{3}:\text{ }\left(
u_{1},u_{2}\right) \longmapsto \left( r_{1}\left( u_{1},u_{2}\right)
,r_{2}\left( u_{1},u_{2}\right) ,r_{3}\left( u_{1},u_{2}\right) \right) , 
\tag{2.3}
\end{equation}%
where $r_{1}$, $r_{2}$ and $r_{3}$ are smooth real-valued functions on a
domain $D\subseteq \mathbb{R}^{2}.$ Denote $g$ the metric on $M^{2}$ induced
from $\mathbb{I}^{3}.$ The components of the first fundamental form of $%
M^{2} $ can be calculated via the induced metric $g$ as follows%
\begin{equation*}
g_{ij}=g\left( r_{u_{i}},r_{u_{j}}\right) ,\text{ }r_{u_{i}}=\frac{\partial r%
}{\partial u_{i}},\text{ }i,j\in \left\{ 1,2\right\} .
\end{equation*}%
The unit normal vector field of $M^{2}$ is completely isotropic, i.e. $%
\left( 0,0,1\right) $. Also, the components of the second fundamental form
are (for details, see \cite{20}, p. 150-155)%
\begin{equation}
t_{ij}=\frac{\det \left( r_{u_{i}u_{j}},r_{u_{1}},r_{u_{2}}\right) }{\sqrt{%
\det \left( g_{ij}\right) }},\text{ }r_{u_{i}u_{j}}=\frac{\partial ^{2}r}{%
\partial u_{i}\partial u_{j}},\text{ }i,j\in \left\{ 1,2\right\} .  \tag{2.4}
\end{equation}%
Thus the \textit{relative curvature} (so-called the \textit{isotropic
curvature} or \textit{isotropic Gaussian curvature}) and the \textit{%
isotropic mean curvature} are respectively defined by%
\begin{equation}
K=\frac{\det \left( t_{ij}\right) }{\det \left( g_{ij}\right) },\text{ }H=%
\frac{g_{11}t_{22}-2g_{12}t_{12}+g_{22}t_{11}}{\det \left( g_{ij}\right) }. 
\tag{2.5}
\end{equation}%
A surface is called \textit{isotropic minimal }(resp. \textit{isotropic flat}%
) if its isotropic mean curvature (resp. relative curvature) vanishes.

In particular, if $M^{2}$ is a Monge surface in $\mathbb{I}^{3}$ of the form 
\begin{equation*}
\left( u_{1},u_{2}\right) \longmapsto \left( u_{1},u_{2},h\left(
u_{1},u_{2}\right) \right) ,
\end{equation*}%
then the metric on $M^{2}$ induced from $\mathbb{I}^{3}$ is given by $%
g_{\ast }=du_{1}^{2}+du_{2}^{2}.$ This implies that $M^{2}$ is always a flat
space with respect to the induced metric $g_{\ast }.$ Thus its Laplacian is
given by%
\begin{equation*}
\bigtriangleup =\sum_{j=1}^{2}\frac{\partial ^{2}}{\partial u_{j}^{2}}.
\end{equation*}%
Also, the matrix of second fundamental form of $M^{2}$ becomes the Hessian
matrix of $h$ (i.e. the square matrix $\left( h_{u_{i}u_{j}}\right) $ of
second-order partial derivatives of the function $h$). Thereby, the formulas
in $\left( 2.5\right) $ reduce to%
\begin{equation}
K=\det \left( h_{u_{i}u_{j}}\right) =h_{u_{1}u_{1}}h_{u_{2}u_{2}}-\left(
h_{u_{1}u_{2}}\right) ^{2},\text{ }H=\bigtriangleup
h=h_{u_{1}u_{1}}+h_{u_{2}u_{2}}.  \tag{2.6}
\end{equation}%
\newpage

\section{Product production functions}

Let us consider the product production function of 2 variables given by 
\begin{equation*}
h:%
\mathbb{R}
_{+}^{2}\longrightarrow 
\mathbb{R}
_{+},\text{ }\left( x,y\right) \longmapsto h\left( x,y\right) =f\left(
x\right) g\left( y\right) ,
\end{equation*}%
where $f$, $g$ are continuous positive real functions with $f^{\prime
}\left( x\right) =\frac{df}{dx}\neq 0$ and $g^{\prime }\left( y\right) =%
\frac{dg}{dy}\neq 0.$ The graph surface $M^{2}$ corresponding to $h$ is of
the form%
\begin{equation}
r\left( x,y\right) =\left( x,y,h\left( x,y\right) =f\left( x\right) g\left(
y\right) \right) ,  \tag{3.1}
\end{equation}%
which we call \textit{product production surface}.

We remark that the surfaces of the form $\left( 3.1\right) $ are known as 
\textit{factorable surfaces }or \textit{homothetical surface}s in ambient
spaces and have been studied in \cite{15,17,23,30}.

The following result provides a complete classification of the product
production surfaces of 2 variables in $\mathbb{I}^{3}$ with constant
relative curvature.

\bigskip

\textbf{Theorem 3.1. }\textit{Let }$M^{2}$\textit{\ be a product production
surface given by }$\left( 3.1\right) $\textit{\ in }$\mathbb{I}^{3}$ \textit{%
with constant relative curvature }$K_{0}$.

\textit{(A) If }$K_{0}=0,$\textit{\ then one of the following occurs:}

\textit{(A.1)\ }$h\left( x,y\right) =c_{1}f\left( x\right) $\textit{\ or }$%
h\left( x,y\right) =c_{2}g\left( y\right) $\textit{\ for nonzero constants }$%
c_{1},c_{2}.$

\textit{(A.2)\ }$h$ \textit{is} \textit{a transcendental production function
of 2 variables} \textit{given by} 
\begin{equation*}
h\left( x,y\right) =A\exp \left( c_{1}x+c_{2}y\right) ,
\end{equation*}%
\textit{\ where }$A,c_{1},c_{2}$\textit{\ are nonzero constants.}

\textit{(A.3) Up to translations of }$x$ \textit{and} $y,$ $h$\textit{\ is a
Cobb-Douglas production function of 2 variables with constant return to
scale.}

\textit{(B) If }$K_{0}\neq 0,$ \textit{then it is negative }$\left(
K_{0}<0\right) $\textit{\ and, up to translations of }$x$ \textit{and} $y,$ $%
h$ \textit{is a Cobb-Douglas production function of 2 variables given by }%
\begin{equation*}
h\left( x,y\right) =\left( -K_{0}\right) xy.
\end{equation*}

\bigskip

\textbf{Proof. }Assume that $M^{2}$ has constant relative curvature $K_{0}$
in $\mathbb{I}^{3}.$ Then, it follows from $\left( 2.6\right) $ that%
\begin{equation}
\left( f^{\prime \prime }g^{\prime \prime }\right) fg-\left( f^{\prime
}g^{\prime }\right) ^{2}=K_{0}  \tag{3.2}
\end{equation}%
where $f^{\prime }=\frac{df}{dx}$ and $g^{\prime }=\frac{dg}{dy}$, etc. We
divide the proof into two cases:

\bigskip

\textbf{Case (i) }$K_{0}=0.$ Then, from $\left( 3.2\right) ,$ both
situations, $f$ or $g$ constants, are solutions for $\left( 3.2\right) $.
This implies the statement (A.1) of the theorem. Now, let us assume that $f$
and $g$ are non-constant functions. Hence, it follows from $\left(
3.2\right) $ that $f$ and $g$ cannot be linear functions. Thus the equation $%
\left( 3.2\right) $ can be rewritten as%
\begin{equation*}
\frac{f^{\prime \prime }f}{\left( f^{\prime }\right) ^{2}}-\frac{\left(
g^{\prime }\right) ^{2}}{g^{\prime \prime }g}=0,
\end{equation*}%
which yields 
\begin{equation}
\frac{ff^{\prime \prime }}{\left( f^{\prime }\right) ^{2}}=\lambda =\frac{%
\left( g^{\prime }\right) ^{2}}{g^{\prime \prime }g}  \tag{3.3}
\end{equation}%
for a nonzero constant $\lambda .$ In order to solve $\left( 3.3\right) $ we
have to distinguish two situations.

\textbf{Case (i.1) }$\lambda =1.$ Then after solving $\left( 2.3\right) $ we
get%
\begin{equation*}
f\left( x\right) =c_{1}\exp \left( c_{2}x\right) \text{ and }g\left(
y\right) =c_{3}\exp \left( c_{4}y\right) ,
\end{equation*}%
where $c_{i}$\ are nonzero constants, $1\leq i\leq 4,$ which gives the the
statement (A.2) of the theorem.

\textbf{Case (i.2) }$\lambda \neq 1.$ Solving $\left( 3.3\right) $ yields 
\begin{equation*}
f\left( x\right) =\left[ \left( 1-\lambda \right) \left( c_{1}x+d_{1}\right) %
\right] ^{\frac{1}{1-\lambda }}\text{ and }g\left( y\right) =\left[ \left( 
\frac{\lambda -1}{\lambda }\right) \left( c_{2}y+d_{2}\right) \right] ^{%
\frac{\lambda }{\lambda -1}}
\end{equation*}%
for nonzero constants $c_{3},c_{4}$ and some constants $d_{1},d_{2}.$ Up to
suitable translations of $x,y$, we obtain%
\begin{equation*}
h\left( x,y\right) =f\left( x\right) g\left( y\right) =Ax^{\frac{1}{%
1-\lambda }}y^{-\frac{\lambda }{1-\lambda }}
\end{equation*}%
for $A=\left[ c_{1}\left( 1-\lambda \right) \right] ^{\frac{1}{1-\lambda }}%
\left[ c_{2}\left( \frac{\lambda -1}{\lambda }\right) \right] ^{\frac{%
\lambda }{\lambda -1}}.$ This proves the statement (A.3) of the theorem.

\bigskip

\textbf{Case (ii) }$K_{0}\neq 0.$ Suppose that\textbf{\ }$f$ and $g$ are
non-linear functions. Hence, we can rewrite $\left( 3.2\right) $ as%
\begin{equation}
\frac{ff^{\prime \prime }}{\left( f^{\prime }\right) ^{2}}-\frac{\left(
g^{\prime }\right) ^{2}}{gg^{\prime \prime }}=\frac{K_{0}}{\left( f^{\prime
}\right) ^{2}gg^{\prime \prime }}.  \tag{3.4}
\end{equation}%
Differentiating of $\left( 3.4\right) $ with respect to $y$ gives%
\begin{equation}
-\left( \frac{\left( g^{\prime }\right) ^{2}}{gg^{\prime \prime }}\right)
^{\prime }=\frac{K_{0}}{\left( f^{\prime }\right) ^{2}}\left( \frac{1}{%
gg^{\prime \prime }}\right) ^{\prime }.  \tag{3.5}
\end{equation}%
From $\left( 3.5\right) ,$ if 
\begin{equation*}
\left( \frac{1}{gg^{\prime \prime }}\right) ^{\prime }=0,
\end{equation*}%
i.e. $gg^{\prime \prime }$ is a nonzero constant $c$ in $\left( 3.5\right) ,$
then we get $\frac{1}{c}\left( \left( g^{\prime }\right) ^{2}\right)
^{\prime }=0,$ which is not possible because $g$ is non-linear function. In $%
\left( 3.5\right) ,$ if 
\begin{equation}
\left( \frac{\left( g^{\prime }\right) ^{2}}{gg^{\prime \prime }}\right)
^{\prime }=0,  \tag{3.6}
\end{equation}%
then $\left( 3.5\right) $ reduces to%
\begin{equation*}
\frac{K_{0}}{\left( f^{\prime }\right) ^{2}}\left( \frac{1}{gg^{\prime
\prime }}\right) ^{\prime }=0,
\end{equation*}%
which yields $gg^{\prime \prime }=d$, $d\neq 0.$ Considering this in $\left(
3.6\right) $ implies $\frac{2}{d}g^{\prime }g^{\prime \prime }=0$ and it is
a contradiction. Thereby we can rewrite $\left( 3.5\right) $ as%
\begin{equation}
-\frac{\left( \frac{\left( g^{\prime }\right) ^{2}}{gg^{\prime \prime }}%
\right) ^{\prime }}{\left( \frac{1}{gg^{\prime \prime }}\right) ^{\prime }}=%
\frac{K_{0}}{\left( f^{\prime }\right) ^{2}}.  \tag{3.7}
\end{equation}%
Since $f$ is a non-linear function, the right-side of $\left( 3.7\right) $
is a function of $x.$ However the left-side of $\left( 3.7\right) $ is
either a constant or a function of $y.$ Both cases are not possible.

Now let either $f$ or $g$ be a linear function. Without loss of generality,
we may assume that $f$ is a linear function, i.e. $f\left( x\right)
=c_{1}x+d_{1},$ $c_{1}\neq 0,$ $d_{1}\in 
\mathbb{R}
$. Then we get from $\left( 3.2\right) $%
\begin{equation*}
g^{\prime }=\frac{\sqrt{-K_{0}}}{c_{1}},\text{ }K_{0}<0.
\end{equation*}%
This implies that $g$ is also a linear function, i.e. $g\left( y\right) =%
\frac{\sqrt{-K_{0}}}{c_{1}}y+d_{2},$ $d_{2}\in 
\mathbb{R}
.$ Thus, up to suitable translations of $x$ and $y,$ we derive%
\begin{equation*}
h\left( x,y\right) =\sqrt{-K_{0}}xy.
\end{equation*}%
This gives of the statement (B) of the theorem.

Therefore, the proof is completed.

\bigskip

Next classifies the product production surfaces of constant isotropic mean
curvature in $\mathbb{I}^{3}.$

\bigskip

\textbf{Theorem 3.2. }\textit{Let }$M^{2}$\textit{\ be a product production
surface given by }$\left( 3.1\right) $\textit{\ in }$\mathbb{I}^{3}$ \textit{%
with constant isotropic mean curvature }$H_{0}$\textit{.\ Then either}

\textit{(A) }%
\begin{equation*}
h\left( x,y\right) =\frac{H_{0}}{g_{0}}x^{2}+d_{1}x+d_{2},\mathit{(or\ }%
\frac{H_{0}}{f_{0}}y^{2}+d_{1}y+d_{2}\mathit{)}
\end{equation*}%
\textit{\ for }$f_{0},g_{0}\in 
\mathbb{R}
-\left\{ 0\right\} ,$ $d_{1},d_{2}\in 
\mathbb{R}
,$ \textit{or}

\textit{(B) }$M^{2}$\textit{\ is isotropic minimal, i.e. }$H_{0}=0,$\textit{%
\ and up translations of }$x$ \textit{and}$\ y,$\textit{\ }$h$ \textit{is a
Cobb-Douglas production function of 2 variables given by }%
\begin{equation*}
h\left( x,y\right) =Axy,A>0.
\end{equation*}

\bigskip

\textbf{Proof. }Assume that $M^{2}$ has constant isotropic mean curvature $%
H_{0}$. Then, by $\left( 2.6\right) ,$ we get%
\begin{equation}
H_{0}=f^{\prime \prime }g+fg^{\prime \prime }.  \tag{3.8}
\end{equation}%
It follows from $\left( 3.8\right) $ that when $g$ is a nonzero constant $%
g_{0}$ we have%
\begin{equation*}
f\left( x\right) =\frac{H_{0}}{g_{0}}x^{2}+d_{1}x+d_{2},\text{ }%
d_{1},d_{2}\in 
\mathbb{R}
,
\end{equation*}%
and analogously if $f$ is a nonzero constant $f_{0},$ we deduce%
\begin{equation*}
g\left( y\right) =\frac{H_{0}}{f_{0}}y^{2}+d_{3}y+d_{4},\text{ }%
d_{3},d_{4}\in 
\mathbb{R}
,
\end{equation*}%
which proves the statement (A) of the theorem.

Now suppose that $f,g$ are non-constant functions. Then $\left( 3.8\right) $
can be rewritten as%
\begin{equation}
\frac{f^{\prime \prime }}{f}+\frac{g^{\prime \prime }}{g}=\frac{H_{0}}{fg}. 
\tag{3.9}
\end{equation}%
After taking the partial derivative of $\left( 3.9\right) $ with respect to $%
x$, we deduce%
\begin{equation}
\left( \frac{f^{\prime \prime }}{f}\right) ^{\prime }\frac{f^{2}}{f^{\prime }%
}=-H_{0}\frac{1}{g}.  \tag{3.10}
\end{equation}%
The left-side of $\left( 3.10\right) $ is etiher a constant or a function of 
$x$ while the other side is a function of $y.$ This case is only possible
when $H_{0}=0$ and 
\begin{equation}
\left( \frac{f^{\prime \prime }}{f}\right) ^{\prime }=0.  \tag{3.11}
\end{equation}%
Similarly, taking the partial derivative of $\left( 3.9\right) $ with
respect to $y,$ we find $H_{0}=0$ and%
\begin{equation}
\left( \frac{g^{\prime \prime }}{g}\right) ^{\prime }=0.  \tag{3.12}
\end{equation}%
This means that $M^{2}$ is isotropic minimal, i.e. $H_{0}=0$. Now let us
assume that $f$ and $g$ are non-linear functions. By solving $\left(
3.11\right) $ and $\left( 3.12\right) ,$ we derive%
\begin{equation}
f\left( x\right) =\beta _{1}e^{\alpha _{1}x}+\beta _{2}e^{-\alpha _{1}x}%
\text{ and }g\left( y\right) =\beta _{3}e^{\alpha _{2}y}+\beta
_{4}e^{-\alpha _{2}y}  \tag{3.13}
\end{equation}%
for nonzero constants $\alpha _{i},\beta _{j},$ $i=1,2$ and $j=1,2,3,4.$
Substituting $\left( 3.13\right) $ into $\left( 3.9\right) $ gives%
\begin{equation*}
\alpha _{1}^{2}+\alpha _{2}^{2}=0,
\end{equation*}%
which is not possible. Then the functions $f$ and $g$ are linear and up to
translations $x$ and $y,$ we obtain $h\left( x,y\right) =Axy,$ $A\neq 0.$
This gives the proof.

\bigskip

\textbf{Remark 3.3. }For the product production function given by $h\left(
x,y\right) =f\left( x\right) g\left( y\right) ,$ we have that $f$ and $g$
are nonconstant functions. Hence, while the statement (A.1) of Theorem 3.1
and the statement (A) of Theorem 3.2 are correct in mathematical
perspective, in reality such product production functions do not exist.

\bigskip

Now let us consider Spillman-Mitscherlich and transcendental production
functions of 2 variables respectively given by 
\begin{equation}
\left\{ 
\begin{array}{l}
h\left( x,y\right) =A\left[ 1-\exp \left( ax\right) \right] \cdot \left[
1-\exp \left( by\right) \right] , \\ 
\left( x,y\right) \in 
\mathbb{R}
_{+}^{n},\text{ }A>0,\text{ }a,b<0%
\end{array}%
\right.  \tag{3.14}
\end{equation}%
and%
\begin{equation}
\left\{ 
\begin{array}{l}
h\left( x,y\right) =Ax^{a_{1}}\exp \left( b_{1}x_{1}\right) \cdot
x_{2}^{a_{2}}\exp \left( b_{2}x_{2}\right) , \\ 
\left( x,y\right) \in 
\mathbb{R}
_{+}^{n},\text{ }A>0,\text{ }a_{i},b_{i}\in 
\mathbb{R}
,\text{ }a_{i}^{2}+b_{i}^{2}\neq 0,\text{ }i=1,2.%
\end{array}%
\right.  \tag{3.15}
\end{equation}

From Theorem 3.1 and Theorem 3.2, we obtain the following results for the
surfaces corresponding to these production functions.

\bigskip

\textbf{Corollary 3.4.}\textit{\ Let }$h$\textit{\ be a
Spillman-Mitscherlich production function of 2 variables given by }$\left(
3.14\right) $\textit{. Then the corresponding graph surface }$M^{2}$\textit{%
\ of }$h$\textit{\ has neither constant relative nor constant isotropic
curvature in} $\mathbb{I}^{3}.$

\bigskip

\textbf{Corollary 3.5. }\textit{Let }$h$\textit{\ be a transcendental
production function of 2 variables given by }$\left( 3.15\right) $ \textit{%
and }$M^{2}$\textit{\ its associated graph surface in }$\mathbb{I}^{3}$%
\textit{. Then:}

\textit{(A) }$M^{2}$ \textit{has constant relative curvature }$K_{0}$ 
\textit{in} $\mathbb{I}^{3}$ \textit{if and only if }$K_{0}=0$\textit{\ and
one of the following occurs:}

\textit{(A.1) }$\alpha _{1}=\alpha _{2}=0$\textit{\ and }$b_{1}\neq 0\neq
b_{2,}$\textit{\ or}

\textit{(A.2) }$\alpha _{1}+\alpha _{2}=1,$\textit{\ }$a_{1}\neq 0\neq a_{2}$
\textit{and }$b_{1}=b_{2}=0.$

\textit{(B) }$M^{2}$ \textit{has constant isotropic curvature }$H_{0}$ 
\textit{in} $\mathbb{I}^{3}$ \textit{if and only if }$H_{0}=0$\textit{\ and }%
$a_{1}=a_{2}=1$ $b_{1}=b_{2}=0.$

\bigskip 

M.E. Aydin

Department of Mathematics

Firat University

23119 Elazig

Turkey

E-mail: meaydin@firat.edu.tr

\bigskip

M. Ergut

Department of Mathematics

Namik Kemal University

59 000 Tekirdag

Turkey

E-mail: mergut@nku.edu.tr


\begin{thebibliography}{99}
\bibitem{1} H. Alodan, B.-Y. Chen, S. Deshmukh, G.E. Vilcu, \textit{On some
geometric properties of quasi-product production models}, arXiv:1512.05190v1
[math.DG], 2015.

\bibitem{2} M.E. Aydin, M. Ergut, \textit{Hessian determinants of composite
functions with applications for production functions in economics},
Kragujevac J. Math. \textbf{38(2)} (2014), 259-268.

\bibitem{3} M.E. Aydin, A. Mihai, \textit{Classifications of quasi-sum
production functions with Allen determinants}, Filomat \textbf{29(6)}
(2015), 1351--1359.

\bibitem{4} M.E. Aydin, M. Ergut, \textit{Composite functions with Allen
determinants and their applications to production models in economics},
Tamkang J. Math. \textbf{45(4)} (2014), 427-435.

\bibitem{5} R. G. Chambers, Applied Production Analysis, Cambridge
University Press, 1998.

\bibitem{6} B.-Y. Chen, G. E. V\^{\i}lcu, \textit{Geometric classifications
of homogeneous production functions}, Appl. Math. Comput. \textbf{225}
(2013), 345--351.

\bibitem{7} B.-Y. Chen, S. Decu, L. Verstraelen, \textit{Notes on isotropic
geometry of production models,} Kragujevac J. Math. \textbf{38(1)} (2014),
23--33.

\bibitem{8} B.-Y. Chen, \textit{On some geometric properties of
h-homogeneous production function in microeconomics},\ Kragujevac J. Math. 
\textbf{35(3)} (2011), 343--357.

\bibitem{9} B.-Y. Chen, \textit{On some geometric properties of quasi-sum
production models}, J. Math. Anal. Appl. \textbf{392} (2012), 192--199.

\bibitem{10} B.-Y. Chen, \textit{Classification of h-homogeneous production
functions with constant elasticity of substitution}, Tamkang J. Math. 
\textbf{43} (2012), 321--328.

\bibitem{11} B.-Y. Chen,\ \textit{Solutions to homogeneous Monge-Ampere
equations of homothetic functions and their applications to production
models in ecenomics}, J. Math. Anal. Appl.\textbf{\ 411} (2014), 223--229.

\bibitem{12} C. W. Cobb, P. H. Douglas, \textit{A theory of production},
Amer. Econom. Rev. \textbf{18 }(1928), 139--165.

\bibitem{13} S. Decu, L. Verstraelen,\textit{\ A note on the isotropical
geometry of production surfaces}, Kragujevac J. Math. \textbf{37(2)} (2013),
217--220.

\bibitem{14} Y. Fu, \textit{Geometric characterizations of quasi-product
production models in economics}, Filomat (2016), to be published.

\bibitem{15} R. Lopez , M. Moruz, \textit{Translation and homothetical
surfaces in Euclidean space with constant curvature}, J. Korean Math. Soc. 
\textbf{52(3)} (2015), 523---535.

\bibitem{16} L. Losonczi, \textit{Production functions having the CES
property}, Acta Math. Acad. Paedagog. Nyh\'{a}i. (N.S.) \textbf{26(1)}
(2010), 113--125.

\bibitem{17} H. Meng, H. Liu, \textit{Factorable surfaces in Minkowski space}%
, Bull. Korean Math. Soc. \textbf{46(1)} (2009), 155--169.

\bibitem{18} S.K. Mishra, \textit{A brief history of production functions},
IUP J. Manage. Econom.\ \textbf{8(4)} (2010), 6--34.

\bibitem{19} H.Pottmann, P. Grohs, N.J. Mitra, \textit{Laguerre minimal
surfaces, isotropic geometry and linear elasticity}, Adv. Comput. Math. 31
(2009), 391--419.

\bibitem{20} H. Sachs,\ Isotrope Geometrie des Raumes, Vieweg, 1990.

\bibitem{21} Z.M. Sipus, \textit{Translation surfaces of constant curvatures
in a simply isotropic space}, Period. Math. Hung. \textbf{68} (2014),
160--175.

\bibitem{22} A. Thompson, Economics of the Firm, Theory and Practice, 3rd
edition, Prentice Hall, 1981.

\bibitem{23} I. Van de Woestyne, \textit{Minimal homothetical hypersurfaces
of a semi-Euclidean space}, Results Math. \textbf{27} (1995), 333--342.

\bibitem{24} G.E. V\^{\i}lcu,\ \textit{A geometric perspective on the
generalized Cobb--Douglas production functions}, Appl. Math. Lett. 24\
(2011), 777--783.

\bibitem{25} A. D. Vilcu, G. E. Vilcu, \textit{On some geometric properties
of the generalized CES production functions}, Appl. Math. Comput. \textbf{218%
} (2011), 124--129.

\bibitem{26} A. D. Vilcu, G. E. Vilcu, \textit{On homogeneous production
functions with proportional marginal rate of substitution}, Math. Probl.
Eng. 2013 (2013), Article ID 732643, 5 pages.

\bibitem{27} A.D. Vilcu, G.E. Vilcu, \textit{Some characterizations of the
quasi-sum production models with proportional marginal rate of substitution}%
, C. R. Math. Acad. Sci. Paris \textbf{353} (2015), 1129-1133.

\bibitem{28} X. Wang, Y. Fu, \textit{Some characterizations of the
Cobb-Douglas and CES production functions in microeconomics}, Abstr. Appl.
Anal. 2013 (2013), Article ID 761832, 6 pages.

\bibitem{29} X. Wang, \textit{A geometric characterization of homogeneous
production models in economics}, Filomat (2016), to be published.

\bibitem{30} Y. Yu, H. Liu, \textit{The factorable minimal surfaces},
Proceedings of The Eleventh International Workshop on Diff. Geom. 11 (2007),
33-39.
\end{thebibliography}
\end{document}